\pdfoutput=1
\documentclass[12pt]{amsart}
\usepackage[utf8]{inputenc}
\usepackage[margin=1 in]{geometry}
\usepackage[sort&compress,square,numbers]{natbib}
\usepackage{graphicx,amssymb,amstext,amsmath,wrapfig,url,bm}
\usepackage[colorlinks]{hyperref}
\usepackage[hang]{subfigure}

\renewcommand{\Pr}[1]{\ensuremath{\text{\rm \textbf{Pr}}\left[{#1}\right]}}
\newcommand{\E}[1]{\ensuremath{\text{\rm \textbf{E}}\left[{#1}\right]}}
\newcommand{\Bin}[2]{\ensuremath{\operatorname{Bin}\left({#1},{#2}\right)}}
\newcommand{\G}[2]{\ensuremath{\mathbb{G}\left({#1},{#2}\right)}}
\renewcommand{\vec}[1]{\mathbf{#1}}

	\def\lemmaD#1{
	\begin{lemma}
	\label{lem:#1}
	}

	\newcommand{\theoremD}[1]{
	\begin{theorem}
	\label{theorem:#1}
	}

	\newcommand{\factD}[1]{
	\begin{fact}
	\label{fact:#1}
	}

	\newcommand{\corD}[1]{
	\begin{cor}
	\label{cor:#1}
	}

	\usepackage{amsmath,amsthm}
	\usepackage{theomac}
	\newcommand{\lemlab}[1]{\label{lemma:#1}}
	\newcommand{\theolab}[1]{\label{theo:#1}}

	\newcommand{\eqlab}[1]{\label{eq:#1}}

	\newcommand{\proplab}[1]{\label{prop:#1}}

	\newcommand{\lemref}[1]{Lemma \ref{lemma:#1}}
	\newcommand{\theoref}[1]{Theorem \ref{theo:#1}}

	\newcommand{\figref}[1]{Figure \ref{fig:#1}}
	\renewcommand{\eqref}[1]{(\ref{eq:#1})}

	\newcommand{\propref}[1]{Proposition \ref{prop:#1}}

	\newtheoremWithMacro{theorem}{Theorem} 
	\newtheorem{lemma}[theorem]{Lemma}
	\newtheoremWithMacro{prop}[theorem]{Proposition}
	\newtheorem{cor}[theorem]{Corollary}

	\newtheorem{question}[theorem]{Question}

\begin{document}

\title{Rigid Components of Random Graphs}
\author{Louis Theran}
\address{Department of Computer Science \\ 
University of Massachusetts, Amherst}
\email{theran@cs.umass.edu}
\urladdr{http://cs.umass.edu/~theran/}
% Primary: Random graphs; Secondary: Graph theory relating to computer science

\begin{abstract}
% We study the size and emergence of the largest rigid component
% in an Erdős-Rényi random graph $\G{n}{p}$ in the case where
% $p=c/n$ for a fixed constant $c>0$.  We show that either $\G{n}{c/n}$
% has only rigid components of size at most $3$ or that the 
% largest rigid component has linear size, almost surely.  
% For $c\le 3$ all the rigid components in $\G{n}{c/n}$ are single edges and 
% triangles, almost surely; for $c>4$ the largest rigid component spans 
% at least $80\%$ of the vertices, almost surely.
The planar rigidity problem asks, given a set of m pairwise distances among a set $P$ of $n$ unknown points, whether it is possible to reconstruct $P$, up to a finite set of possibilities (modulo rigid motions of the plane).  The celebrated Maxwell-Laman Theorem from Rigidity Theory says that, generically, the rigidity problem has a combinatorial answer: the underlying combinatorial structure must contain a spanning minimally-rigid graph (Laman graph).  In the case where the system is not rigid, its inclusion-wise maximal rigid substructures (rigid components) are also combinatorially characterized via the Maxwell-Laman theorem, and may be found efficiently.

Physicists have used planar combinatorial rigidity to study the phase transition between liquid and solid in network glasses.  The approach has been to generate a graph via a stochastic process and then to analyze experimentally  its rigidity properties.  Of particular interest is the size of the largest rigid components.

In this paper, we study the emergence of rigid components in an Erdős-Rényi random graph $\G{n}{p}$, using the parameterization $p=c/n$ for a fixed constant $c>0$.  Our first result is that for all $c>0$, almost surely all rigid components have size $2$, $3$ or $\Omega(n)$.  We also show that for $c>4$, almost surely the largest rigid components have size at least $n/10$.  

While the $\G{n}{p}$ model is simpler than those appearing in the physics literature, these results are the first of this type where the distribution is over all graphs on $n$ vertices and the expected number of edges is $O(n)$.
\end{abstract}

\maketitle

%%%%%%%%%%%%%%%%%%%%%%%%%%%%%%%%%%%%%%%%%%%%%%%%%%%%%%%%%%%%%%%%%%%%%
%%%%%%%%%%%%%%%%%%%%%%%%%%%%%%%%%%%%%%%%%%%%%%%%%%%%%%%%%%%%%%%%%%%%%
\section{Introduction}
The problem of the phase transition between liquid and solid states of 
glasses is an important open problem in material physics \cite{anderson:glass:1995}.  Glasses are highly disordered solids that undergo a rapid transition as they 
cool.  
%This as-yet unsolved problem was described by Anderson
%as ``The deepest and most interesting unsolved problem in solid state theory.''

To study the phase transition, Thorpe \cite{JaTh95} proposed a \emph{geometric}
model for the glass problem, in which bonds between the atoms are viewed as  
\textbf{fixed-length bars} (the bonds) connected by \textbf{universal joints} 
(the atoms) with full rotational degrees of freedom.  
Such a structure is called a \textbf{planar bar-and-joint framework} 
(shortly bar-joint framework, or simply framework), and these are fundamental 
objects of study in the field of \textbf{combinatorial rigidity} (see,
e.g., \cite{graver:servatius:rigidityBook:1993} for a survey).

A bar-joint framework is 
\textbf{rigid} if the only continuous motions of the joints
preserving the lengths and connectivity of the bars are rigid 
motions of the plane, and otherwise it is \textbf{flexible}.
When a framework is flexible, it decomposes uniquely 
into inclusion-wise maximal rigid substructures 
which are called \textbf{rigid components} (shortly components);
a component is non-trivial if it is larger than a single edge.
In the planar case, the celebrated Maxwell-Laman Theorem \cite{laman}
gives a complete characterization of \emph{generically} minimally 
rigid bar-joint frameworks in terms of a combinatorial condition, 
which allows rigidity properties to be studied in terms of efficiently
checkable graph properties.

The sequence of papers
 \cite{JaTh95,jacobs:hendrickson:PebbleGame:1997a,thorpe:rigidity:glasses:2002,chubynsky2002rigid,ThJaChRa99} 
studies the emergence of \emph{large} rigid subgraphs in graphs generated
by various stochastic processes, with the edge probabilities and 
underlying topologies used to model the temperature and chemical composition 
of the system.  Two important observations are that: (1) very large 
rigid substructures emerge very rapidly; (2) the transition appears to 
occur slightly below average degree $4$ in the the planar bar-joint model.

\subsection*{Main result novelty}
In this paper, we study the emergence of rigid 
components in random graphs generated by a 
simple, well-known stochastic process: the Erdős-Rényi random graph model $\G{n}{p}$, 
in which each edge is included with probability $p$, independently.  
We consider edge probabilities of the form $p=c/n$, where $c$ is a
fixed constant, and consider the size of the largest 
rigid components in $\G{n}{p}$.  

% Define the random variable $r(n,p)$ to be the size of the 
% largest non-trivial rigid component 
% of $\G{n}{p}$ (this may not be unique).  Our 
% main result is the following statement about $r(n,p)$.
Our main result is the following statement about rigid components in 
$\G{n}{c/n}$.

\begin{theorem}[\emergence][{\bf Size and emergence of a large rigid component}]
	\theolab{emergence}
Let $c>0$ be a constant.  Almost surely, all rigid components in 
$\G{n}{c/n}$ span $2$, $3$, or $\Omega(n)$ vertices. If $c>4$, then 
almost surely there are components of size at least $n/10$.
\end{theorem}
(A random graph has a property almost surely if the 
probability of $\G{n}{p}$ having it tends to one as $n\to \infty$.)

To the best of our knowledge, this is the first proven result on
the emergence of rigid components in random graphs that have,
almost surely, close to $2n-3$ edges (the number required for minimal rigidity)
but \emph{no other special assumptions}, such as being highly connected or a 
subgraph of a hexagonal lattice, both of which play critical roles 
in the previous results on the rigidity of random graphs.  

It is important to note that rigidity is inherently a non-local 
phenomenon: adding a single edge to a graph that has no non-trivial 
rigid components may rigidify the entire graph (or removing a single
edge may cause a large rigid component to shatter).  It is this property
of rigidity that distinguishes it from the well-studied
$k$-core problem in random graph theory.

In the proof of \theoref{emergence}, we formalize the experimental observation
that rigid components, once they appear, are very likely to grow rapidly.  Although
the proof of \theoref{emergence} relies mainly on standard tools for bounding 
sums of independent random variables, our result seems to be the first that 
directly analyzes rigidity properties of $\G{n}{p}$, rather than reducing to
a connectivity property.

%%% put something here

\subsection*{Related work.}
Jackson, et al. \cite{servatius:random:2008} studied the space of random 
$4$-regular graphs and showed that they are almost surely globally rigid 
(see \cite{jackson:jordan:connectedRigMatroids:2005,connelly:global:2005}).
They also established a threshold for $\G{n}{p}$ to be rigid at 
$p=n^{-1}(\log n+2\log\log n+\omega(1))$, which coincides with the threshold 
for $\G{n}{p}$ to almost surely have all vertices with degree at least $2$.
The approach in \cite{servatius:random:2008} is based on combining results on the 
connectivity of random graphs  (e.g., \cite[Theorem 4]{luczak:1991:k-core})
and theorems linking rigidity and connectivity proved in \cite{servatius:random:2008}
and also \cite{lovasz:yemini,jackson:jordan:connectedRigMatroids:2005}.  In the 
$\G{n}{p}$ model, the techniques there seem to rely on the existence of a very large 
$6$-core, so it does not seem that they can be easily adapted to our setting 
when $c$ is close to $4$ (below the threshold for even the $4$-core to emerge 
\cite{spencer:k-core:1996}).

Holroyd \cite{holroyd:98} extended the formal study of connectivity 
percolation \cite{bollobas:percolation:2006} to rigidity percolation 
in the hexagonal lattice.  He shows, via a reduction to 
connectivity percolation, that there is an edge-probability 
threshold for the existence of an infinite\footnote{Rigidity of infinite 
frameworks is a subtle concept, and \cite{holroyd:98} devotes 
careful attention to its development.} 
rigid component in the hexagonal lattice which is higher than 
that for connectivity.  It is also shown in \cite{holroyd:98}
that the infinite component, when it exists, is unique for all 
but a countable set of edge probabilities $p$.  All the 
proofs in \cite{holroyd:98} rely in an essential way on the 
structure of the hexagonal lattice (in particular that a suitably 
defined tree in its dual graph is a dual of a rigid component).

The fundamental $k$-core problem in random graph theory has been 
studied extensively, with a number of complete solutions.  {\L}uczak
\cite{luczak:1991:k-core} first proved that for $k\ge 3$, the (it is always unique, if present) $k$-core
is, almost surely, either empty or has linear size.  
Pittel, et al. 
solved the $k$-core problem, giving an exact threshold for its emergence and 
bounds on its size \cite{spencer:k-core:1996}.  Janson and Luczak gave an alternative 
proof of this result, using simpler stochastic processes \cite{janson:k-core:2007}.  All these 
results are based on analyzing a process that removes low-degree vertices one
at a time, which does not apply in the rigidity setting.

\section{Preliminaries}
In this section we give the technical preliminaries required for the proof
of \theoref{emergence}.

\subsection*{Combinatorial rigidity}
An \textbf{abstract bar-and-joint framework} $(G,\bm{\ell})$ 
is a graph $G=(V,E)$ and vector of non-negative \textbf{edge lengths}
$\bm{\ell}=\ell_{ij}$, for each edge $ij\in E$.  A realization $G(\vec p)$
of the abstract framework $(G,\bm{\ell})$ is an embedding of $G$ onto the 
planar point set $\vec p=(\vec p_i)_{1}^n$ with the property that 
for all edges $ij\in E$, $||\vec p_i-\vec p_j||=\ell_{ij}$.  
The framework $(G,\bm{\ell})$ is \textbf{rigid} if it has only a discrete set of 
realizations modulo trivial plane motions, and is \textbf{flexible} otherwise.

A graph $G=(V,E)$ is \textbf{$(2,3)$-sparse} 
if every subgraph induced by 
$n'\ge 2$ vertices has at most $2n'-3$ edges.  If, in addition, 
$G$ has $2n-3$ edges, $G$ is \textbf{$(2,3)$-tight} (shortly, Laman).

The Maxwell-Laman Theorem completely characterizes the rigidity of 
generic planar bar-joint frameworks.
\begin{prop}[\laman][\textbf{Maxwell-Laman Theorem \cite{laman}}]
A generic bar-joint framework in the plane is minimally 
rigid if and only if its graph is $(2,3)$-tight.
\end{prop}
Genericity is a subtle concept, and we refer the reader to our
paper \cite{streinu:theran:LamanGenericity:2008} for a detailed 
discussion.  In the following it suffices to note that for a fixed $G$ 
almost all $\vec p$ are generic, and that, by the Maxwell-Laman Theorem, 
all generic frameworks $G(\vec p)$ have the same rigidity properties.

\begin{figure}[htbp]%{l}{0.4\textwidth}
    \centering
%\vspace{-0.2 in}
    \subfigure[]{\includegraphics[width=0.4\textwidth]{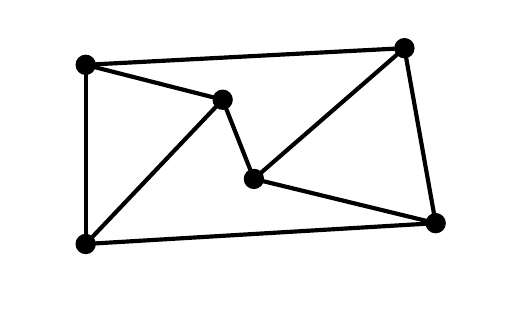}}%\vspace{-0.2 in}
    \subfigure[]{\includegraphics[width=0.4\textwidth]{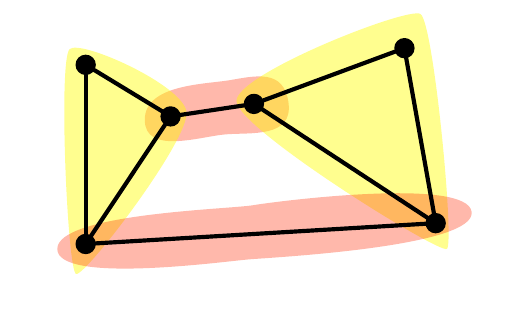}}
    \caption{Laman graphs and rigid components: (a) a Laman graph on $n=6$
vertices; (b) a flexible graph with its rigid components indicated.}
    \label{fig:laman-examples}
\end{figure}

If $G$ contains a spanning Laman graph it is \textbf{$(2,3)$-spanning} 
(shortly rigid).  
A rigid induced subgraph is called a \textbf{spanning block}
(shortly block), and 
an inclusion-wise maximal block is a \textbf{spanning component} 
(shortly component)\footnote{In \cite{streinu:lee:pebbleGames:2008}
the terms ``block'' and ``component'' are
reserved for induced subgraphs of Laman graphs, but there is 
no concern of confusion here.}.  By  
\cite[Theorem 5]{streinu:lee:pebbleGames:2008}, every graph decomposes uniquely into 
components, and every edge is spanned by exactly one component.  A 
component is non-trivial if it contains more than one edge. \figref{laman-examples}(a)
shows and example of a Laman graphs.  \figref{laman-examples}(b) has an example
of a flexible graph with its components indicated; they are the two triangles and two 
trivial components consisting of a single edge only.

An alternative characterization of Laman graphs is via so-called \textbf{Henneberg constructions}, which are local moves that transform Laman graphs on $n$
vertices to Laman graphs on $n+1$ vertices (see \cite[Section 6]{streinu:lee:pebbleGames:2008}).  The \textbf{Henneberg I} move adds a new
vertex $n$ to a Laman graph $G$ and attaches it to two neighbors in $V(G)-n$.
It is a fundamental result of rigidity theory that the Henneberg I move 
preserves generic rigidity \cite{laman}\footnote{This fact, along with 
an analogous result for the so-called Henneberg II move, which adds a vertex
of degree 3, is the core of Laman's proof.}.

We summarize the properties of rigid graphs and components that we will use 
below in the following proposition.
\begin{prop}[\xyzzy][{\bf Properties of rigid graphs and rigid components}]
	\proplab{rigidprops}
	Let $G=(V,E)$ be a simple graph with $n$ vertices.
	\begin{itemize}
		\item[(a)] $G$ decomposes uniquely into rigid components (inclusion-wise 
		maximal induced Laman graphs), and every edge is in some component \cite[Theorem 5]{streinu:lee:pebbleGames:2008}.
		\item[(b)] Adding an edge to a graph $G$ never decreases the size of any 
		rigid component \cite[Theorem 2]{streinu:lee:pebbleGames:2008}.
		\item[(c)] If $G'$ is a block in $G$ with vertices $V'\subset V$ 
		and there is a vertex $i\notin V'$ with at least two neighbors in 
		$V'$, then $G'$ is not a component of $G$.
		\item[(d)] If $G$ has at least $2n-2$ edges, then it contains a
		component spanning at least $4$ vertices \cite[Theorem 2 and Theorem 5]{streinu:lee:pebbleGames:2008}.
	\end{itemize}
\end{prop}

What we have presented here is a small part of a well-developed combinatorial
and algorithmic theory of $(k,\ell)$-sparse graphs.  We refer the reader to 
\cite{streinu:lee:pebbleGames:2008} for a detailed treatment of the 
rich properties of sparse graphs.

\subsection*{Tools from random graph theory}
One of our main technical tools is the following result
on the size of dense subgraphs in $\G{n}{c/n}$  
due to {\L}uczak \cite{luczak:1991:k-core}.  Since it appears 
without proof in \cite{luczak:1991:k-core}, we give our own 
in the appendix.
\begin{prop}[\density][{\bf Density Lemma \cite{luczak:1991:k-core}}]
	\proplab{density}
Let $a$ and $c$ be real constants with $a>1$ and $c>a$.  
Almost surely, $\mathbb{G}(n,c/n)$ has no subgraphs with 
at most $k=t(a,c)n$ vertices and at least $akn$ edges, where
\[
t(a,c)=\left(\frac{2a}{c}\right)^{\frac{a}{a-1}}
   e^{-\frac{a+1}{a-1}}
\]
\end{prop}

We will also make use of a fairly general form of the Chernoff bound 
for the upper tail of the binomial.
\begin{prop}[\chernoff][Chernoff bound]\proplab{chernoff}
	Let $\Bin{N}{p}$ be a binomial random variable with parameters $n$
	and $p$.  Then for all $\delta>0$,
	\[
		\Pr{\Bin{N}{p}\ge (1+\delta)Np} \le 
			\left(\frac{e^\delta}{(1+\delta)^{(1+\delta)}}\right)^{Np}
	\]
\end{prop}
Large deviation bounds of this type are attributed to Chernoff
 \cite{chernoff:tails:1952}, and are standard in combinatorics.  The 
specific form of \propref{chernoff} appears in, e.g., 
\cite[Theorem 4.1, p. 68]{motwani:randombook:1995}.

\section{Proofs}
In this section we prove the main result of this paper.  

\emergence

\subsection*{Proof outline.}
Here is the proof strategy in a nutshell.  Because any rigid component
with $n'\ge 4$ vertices
must be somewhat dense, the very general bound of \propref{density}
implies that for $p=c/n$ all the components are either trivial, triangles, 
or spanning a constant fraction of the vertices in $\G{n}{c/n}$ (\lemref{linearsize}).  We then improve upon our bounds on the probability of components of size
$sn$, for $s\in (0,1)$ by formalizing the observation that 
such components are likely to ``grow'' (\lemref{compprob}) and then 
optimizing $s$ (\lemref{lowerbound}).

The rest of this section contains the details.

\subsection*{Rigid components have either constant or linear size}

We start by proving that non-trivial rigid components are all
very large or triangles, almost surely. 
\begin{lemma}\lemlab{linearsize}
Let $c>0$ be a fixed constant.  Almost surely, all rigid components in 
$\G{n}{c/n}$ have size $2$, $3$, or $\Omega(n)$.
\end{lemma}
\begin{proof}
By \propref{rigidprops}(a), any rigid component on $n'\ge 4$ vertices
has at least $\frac{5}{4}n'$ edges (with equality for $n'=4$).  The 
lemma then follows from \propref{density} and the well-known fact that
almost surely $\G{n}{c/n}$ contains a triangle \cite[Theorem 4.1, p. 79]{bollobas:randomgraph:book}.
\end{proof}
{\bf Remark:} In fact, this proof via \propref{density} implies a stronger result, which is that almost surely $\G{n}{c/n}$ does not contain any sub-linear size induced subgraphs
with enough edges to be non-trivial rigid blocks, except for triangles.  

For $c>4$, the number of edges in $\G{n}{c/n}$ implies that it has at least
one large rigid component, almost surely.
\begin{lemma}\lemlab{cge4}
Let $c>4$.  Almost surely, $\G{n}{c/n}$ contains at least one component 
of size $\Omega(n)$.
\end{lemma}
\begin{proof}
For any $\epsilon>0$ $\G{n}{(4+\epsilon)/n}$ has at least $2n-2$ edges with high 
probability.
\propref{rigidprops}(d) then implies that almost surely $\G{n}{(4+\epsilon)/n}$
contains at least one rigid component with at least $4$ vertices.  By \lemref{linearsize}, all of these span at least $t(a,4+\epsilon)n$ vertices.

By \propref{rigidprops}(b) the size of rigid components 
is an increasing property and \cite[Theorem 2.1, p. 36]{bollobas:randomgraph:book},
this lower bound on size holds, almost surely, for any $c>4$.
\end{proof}

\subsection*{For $c>4$ the largest component is very large}
We now turn to improving the lower bound on the size of rigid components.  To
do this, we will use the maximality of components as well as their edge density.

\begin{lemma}\lemlab{compprob}
The probability that a fixed set of $k$ vertices spans a component in $\G{n}{c/n}$ is at most
\begin{eqnarray}\eqlab{compprob}
	\Pr{\Bin{k^2/2}{c/n}\ge 2k-3}\left((1-c/n)^k+k\frac{c}{n}(1-c/n)^{k-1}\right)^{n-k}
\end{eqnarray}
\end{lemma}
\begin{proof}
To induce a component, a set $V'$ of $k$ vertices must span at least
$2k-3$ edges by \propref{rigidprops}(a).  By \propref{rigidprops}(c)
if $V'$ spans a component, no vertex outside of $V'$ can have 
more than one neighbor in $V'$.  The two terms in \eqref{compprob}
correspond to these two events, which are independent.
\end{proof}
{\bf Remark:} This estimate of the probability of a set of vertices inducing a 
component is very weak, since it uses only the number of edges induced by $V'$ (not
their distribution)
and the simplest local obstacle to maximality.  Any improvement in this part of 
the argument would translate into improvements in the lower bound on the 
size of components.

\begin{lemma}\lemlab{lowerbound}
For $c>4$, almost surely  all components span at least $n/10$ vertices.
\end{lemma}
\begin{proof}
With the assumptions of the lemma, by \lemref{cge4}, $\G{n}{c/n}$ almost surely 
has no blocks of size smaller than $tn$, where $t$ is a constant
independent of $n$.  It follows from 
\propref{rigidprops}(a) that $\G{n}{c/n}$ almost surely has no components
smaller than $tn$.

Let $X_k$ to be the number of components of size $k$ and let $s$ be a parameter
to be selected later.  We will show that $\sum_{k=4}^{sn}\E{X_k}=o(1)$,
which implies the lemma by a Markov's inequality.  As noted above, 
$\sum_{k=4}^{tn}\E{X_k}=o(1)$, so we concentrate on $k\in [tn, sn]$.

By \lemref{compprob}
\[
		\begin{split}
			\E{X_k} & \le \binom{n}{k}\Pr{\Bin{k^2/2}{c/n}\ge 2k-3}\left((1-c/n)^k+k\frac{c}{n}(1-c/n)^{k-1}\right)^{n-k} \\
		& \le \left(\frac{en}{k}\right)^k \Pr{\Bin{k^2/2}{c/n}\ge 2k}
		\left((1-c/n)^k+k\frac{c}{n}(1-c/n)^{k-1}\right)^{n-k} +o(1) \\
		\end{split}
\]

Setting $k=sn$ and letting $c=4+\epsilon$, we use the Chernoff bound 
to obtain
\begin{multline*}
		\left(\frac{e}{s}\right)^{sn} \Pr{\Bin{k^2/2}{(4+\epsilon)/n}\ge 2sn}
		\left((1-(4+\epsilon)/n)^sn +cs(1-(4+\epsilon)/n)^{sn-1}\right)^{n-sn} \le\\
		\left(\frac{e}{s}\right)^{sn}
		\left(e^{\frac{-\epsilon  s-4 s+4}{s (\epsilon +4)}}
		   \left(\frac{-\epsilon  s-4 s+4}{s (\epsilon
		   +4)}+1\right)^{-\frac{-\epsilon  s-4 s+4}{s (\epsilon
		   +4)}-1}\right)^{\frac{1}{2} n s^2 (\epsilon +4)}
		\left(e^{-(4+\epsilon)s}(1+(4+\epsilon)s)\right)^{n-sn}
\end{multline*}
As $\epsilon\to 0$ the right-hand side approaches
\[
	e^{n s} \left(e^{\frac{1}{s}-1} \left(\frac{1}{s}\right)^{-1/s}\right)^{2
	   n s^2} \left(\frac{1}{s}\right)^{n s} \left(e^{-4 s} (4
	   s+1)\right)^{n-n s}
\]
Substituting $s=1/10$, this simplifies to 
\[
	2^{-n/10} 5^{-n} 7^{9 n/10} e^{-2 n/25} = e^{-\Theta(n)}
\]
(which can be seen by taking the logarithm and factoring out $n$).  
Since this 
bound is good for any $s'\in [t,{1/10}]$, we have 
$\sum_{k=tn}^{n/10} \E{X_k}\le n e^{-\Theta(n)} = o(1)$.

By \propref{rigidprops}(b) the size of rigid components 
is an increasing property and \cite[Theorem 2.1, p. 36]{bollobas:randomgraph:book},
this lower bound on size holds almost surely for any $c>4$.
\end{proof}

\section{Conclusions and open problems}
We considered the question of the size and emergence of rigid components in a random 
graph $\G{n}{c/n}$ as $c$ increases, and we proved that almost surely 
all rigid components in $\G{n}{c/n}$ are single edges,  triangles or span $\Omega(n)$ vertices.  
For $c>4$, we proved that, almost surely, the largest rigid components 
span at least $n/10$ vertices.

The most natural open question is whether there is a threshold constant 
for rigid components in $\G{n}{p}$.
\begin{question}[\textbf{Existence of a threshold constant}] 
	Is there a constant $c_r$ at which a linear-sized rigid 
	component appears in $\G{n}{(c_r+\epsilon)/n}$ almost surely, 
	and $\G{n}{(c_r-\epsilon)/n}$ almost surely has no large 
	rigid components?
\end{question}

The other important question is about the structure of large 
rigid components when they emerge.
\begin{question}[\textbf{Structure of large rigid components in $\G{n}{c/n}$}]
	Is there almost surely only one large rigid component in $\G{n}{c/n}$,
	and what are the precise bounds on its size?
\end{question}
We have observed in computer simulations that when linear sized rigid components
are present, there is only one, and it is much larger than 
$n/10$.

%%%%%%%%%%%%%%%%%%%%%%%%%%%%%%%%%%%%%%%%%%%%%%%%%%%%%%%%%%%%%%%%%
% includes all the bibliography files
%

%\input{all.tex}

\appendix

\section{Proof of \propref{density}}
We give a proof of \propref{density}, which specializes 
a technical lemma from \cite{luczak:1991:k-core}.  It appeared there 
without proof.
\density

\begin{proof}
Let $t=t(n)$ be a parameter to be picked later, and say that a subgraph is bad
if it has at most $tn$ vertices and at least $atn$ vertices, and let $X$
be the number of bad subgraphs in $\G{n}{c/n}$.  The proof is via a 
first moment argument.  We observe that for any set of $k$ vertices 
in $\G{n}{c/n}$, the number of induced edges is a 
random variable that is dominated by the binomial 
random variable $\Bin{k^2/2}{c/n}$.  

Let $X_k$ be the number of 
bad subgraphs of size $k$.  By definition,
\begin{eqnarray}\eqlab{boundingsum}
	\E{X}\le \sum_{k=2}^{tn} \E{X_k}\le \sum_{k=2}^{tn} \binom{n}{k}
	 	\sum_{j=ak}^{{k^2}/2} \binom{k}{j}
		\left(\frac{c}{n} \right)^j \left( 1-\frac{c}{n}\right)^{k-j}
\end{eqnarray}
and we will show that the right hand side $o(1)$, for a choice of 
$t$ independent of $n$, which implies the lemma 
since Markov's inequality shows that $\Pr{X>0}\le \E{X}$.
To do this, we split the sum in \eqref{boundingsum} into two parts:
$2\le k\le n^{\epsilon}$, where $\epsilon<\min\{\frac{2a(1-1/a)}{a+3},1/2\}$; 
and $n^\epsilon<k\le tn$.

For the small terms, we start by expanding $\E{X_k}$ directly:
\[
	\begin{split}
		\sum_{k=2}^{n^\epsilon} \binom{n}{k}
		 	\sum_{j=ak}^{{k^2}/2} \binom{{k^2}/2}{j}
			\left(\frac{c}{n}\right)^j\left(1-\frac{c}{n}\right)^{k-j} & \le
		\sum_{k=2}^{n^\epsilon} \left(\frac{en}{k}\right)^k
		 	\sum_{j=ak}^{{k^2}/2} \left(\frac{ek}{2j}\right)^j
			\left(\frac{c}{n}\right)^j
			\\
		& \le 	
			\sum_{k=2}^{n^\epsilon} \left(\frac{en}{k}\right)^k
			\sum_{j=ak}^{{k^2}/2} \left(\frac{ekc}{2an}\right)^{ak} \\
		& \le 	
			\sum_{k=2}^{n^\epsilon} \sum_{j=ak}^{{k^2}/2}
				\left(\frac{k^{ak}}{k^k}\right)
				\left(
					\frac{e^{1+1/a} c}{2 a n^{1-1/a}}
				\right)^{2a} \\
		& \le
			n^{\epsilon(3+a)}
			 	\left(
					\frac{e^{1+1/a} c}{a n^{1-1/a}}
				\right)^{2a} \\
		& = o(1) 
	\end{split}
\]

For $k>n^{\epsilon}$ we parameterize $k$ as $tn$ with $t>n^{\epsilon-1}$
and use the 
use the Chernoff inequality 
to bound the probability that  
$\text{Bin}(\frac{1}{2}(tn)^2,c/n)>atn$.
Plugging in to \propref{chernoff} with $\delta=\frac{2 a}{c t}-1$
shows that the probability of
any particular set of $tn$ vertices inducing a bad subgraph is 
at most
\begin{equation}\eqlab{prdense}
\left( e^{\frac{2 a}{c t}-1} \left(\frac{2a}{c
   t}\right)^{-\frac{2 a}{c t}}\right)^{\frac{1}{2} c n t^2}
\end{equation}

The number of sets of size $tn$ is at most $(e/t)^{tn}$.  Multiplying  
with \eqref{prdense} gives a bound on $\E{X_{tn}}$:
\[
\E{X_{tn}}\le \left( e^{(a+1) t-\frac{c t^2}{2}} \left(\frac{1}{t}\right)^t
   \left(\frac{2 a}{c t}\right)^{-a t}\right)^n
\]
We can show that $t$ can be chosen independently 
of $n$ to make the inner expression strictly less than one.  
Taking the logarithm, we obtain 
\[
	t 
	  \left(
		a + 1 -\frac{ct}{2}-a\log\left(\frac{2a}{c}\right)+(a-1)\log(t)
	  \right)
\]
Plugging in $t(a,c)$ from the statement, this simplifies to 
$-\frac{c t(a,c)^2}{2}<0$, from our assumptions on $a$ and $c$,
and this function decreases with $t$ in the interval
 $[n^{\epsilon-1},t(a,c)]$.

It follows that 
\[
   \sum_{k=n^{\epsilon}}^{t(a,c)n} \E{X_k} \le n e^{-\Theta(n)} = o(1)
\]
completing the proof.
\end{proof}
\end{document}